\documentclass[11pt]{article}

\usepackage{amsmath,amsthm,amsfonts,amssymb,amscd}
\usepackage{dsfont}
\usepackage{subfig}
\usepackage{graphicx}
\usepackage{enumerate}
\usepackage{hyperref}
\usepackage[latin1]{inputenc}

\newtheorem{teo}{Theorem}
\numberwithin{teo}{section}
\newtheorem{prop}[teo]{Proposition}
\newtheorem{lema}[teo]{Lemma}

\newenvironment{pf}{\noindent\textrm{Proof.}}{\hfill $\square$}

\newcommand{\R}{\mathbb{R}}

\usepackage{authblk}
\usepackage{footnote}

\providecommand{\keywords}[1]{\textit{Keywords:} #1}

\date{}

\title{\textbf{On the vertex degrees of the skeleton of the matching polytope of a graph}}

\author[a]{Nair Abreu \thanks{Nair Abreu: nairabreunovoa@gmail.com. The work of this author is partially supported by CNPq, PQ-304177/2013-0 and Universal Project-442241/2014-3.} }
\author[b]{Liliana Costa \thanks{Liliana Costa: lmgccosta@gmail.com. The work of this author is partially supported by the Portuguese Foundation for Science and Technology (FCT) by CIDMA -project UID/MAT/04106/2013.} }
\author[c]{Carlos Henrique do Nascimento \thanks{Carlos Henrique do Nascimento: carloshenrique@id.uff.br.} }
\author[d]{Laura Patuzzi \thanks{Laura Patuzzi: laurapatuzzi@gmail.com.} }

\affil[a]{Engenharia de Produção,  COPPE, Universidade Federal do Rio de Janeiro, Brasil}
\affil[b]{Col\'{e}gio Pedro II, Rio de Janeiro, Brasil}
\affil[c]{Instituto de Ci\^{e}ncias Exatas, Universidade Federal Fluminense, Volta Redonda, Brasil}
\affil[d]{Instituto de Matem\'{a}tica, Universidade Federal do Rio de Janeiro, Brasil}

\begin{document}

\maketitle 

\begin{abstract}
 The convex hull of the set of the incidence vectors of the matchings of a graph $G$ is the matching polytope of the graph,  $\mathcal{M}(G)$. The graph whose vertices and edges are the vertices and edges of $\mathcal{M}(G)$ is the skeleton of the matching polytope of $G$, denoted  $\mathcal{G}(\mathcal{M}(G))$. Since the number of vertices of $\mathcal{G}(\mathcal{M}(G))$ is huge, the structural properties of these graphs have been studied in particular classes. In this paper, for an arbitrary graph $G$, we obtain a formulae to compute the degree of a vertex of $\mathcal{G}(\mathcal{M}(G))$ and prove that the minimum degree of $\mathcal{G}(\mathcal{M}(G))$ is equal to the number of edges of $G$. Also, we identify the vertices of the skeleton with the minimum degree and characterize regular skeletons of the matching polytopes.
\end{abstract}

\keywords{\texttt{graph, matching polytope, degree of matching.}}

\section{Introduction}
Let $G:= G(V(G), E(G))$ be a simple and connected graph with vertex set $V:=V(G)=\{v_1, v_2, \ldots,$ $ v_n\}$ and set of edges $E:=E(G)= \{e_1, e_2, \ldots , e_{m}\}$. For each $k,$ $1 \leq k \leq m$, $e_k = v_iv_j$ is an edge incident to the adjacent vertices $v_i$ and $v_j,$ $1 \leq i < j \leq n$. The set of adjacent vertices of $v_i$ is $N_{G}(v_i)$, called the \textit{neighborhood} of $v_i,$ whose cardinality $d(v_i)$ is the {\it degree} of $v_i$.

If two edges have a  common vertex they are said to be \emph{adjacent edges}. For a given edge $e_k$, $I(e_k)$ denotes the set of adjacent edges of $e_k$. Two non adjacent edges are \emph{disjoint} and  a set of pairwise disjoint edges $M$  is a \emph{matching} of $G$.  An unitary edge set is an \emph{one-edge matching} and the empty set is the \emph{empty matching}, $\varnothing$. A vertex  $v \in V(G)$ is said to be \emph{$M$-saturated} if there is an edge of $M$ incident to $v$. Otherwise, $v$ is said an \emph{$M$-unsaturated} vertex. A \emph{perfect matching} $M$ is one for which every vertex of $G$ is an $M$-saturated.

For a natural number $k$, a \textit{path} with length $k$,  $P_{k+1}$ (or simply $P$), is a sequence of distinct vertices $v_1v_2\dots v_k v_{k+1}$ such that, for $1 \leq i \leq k$, $e_i=v_iv_{i+1}$ is an edge of $G$.  A \textit{cycle} $C=C_k$, with length $k$, is obtained by path $P_k$ adding the edge $v_1v_k$. If $k$ is odd, $C_k$ is said to be an \textit{odd cycle}. Otherwise, $C_k$ is an \textit{even cycle}. When it is clear, we denote a path and a cycle by a sequence of their respective edges $e_1e_2 \dots e_k$ instead of their respective sequences of vertices. Given a matching $M$ in $G$, a path $P$ (or, cycle $C$) is \emph{$M$-alternating path (or, $M$-alternating cycle)} in $G$  if given two adjacent edges of P (or, C), one belongs to M and one belongs to $E\backslash M$. Naturally, the set of edges of $P\backslash M$ (and $C\backslash M$) is also a matching in $G$. For more basic definitions and notations of graphs, see \cite{Bolob, Diestel} and, of matchings, see \cite{Lovasz}.

A \emph{polytope} of $ \R^n$ is the convex hull $\mathcal{P} = conv{\{x_1,\ldots,x_r\}}$ of a finite set of vectors $x_1,\ldots,x_r \in \R^n$. Given a polytope $\mathcal{P}$, the \emph{skeleton} of $\mathcal{P}$ is a graph $\mathcal{G}(\mathcal{P})$ whose vertices and edges are, respectively, vertices and edges of $\mathcal{P}$.

Ordering the set $E$ of $m$ edges of $G$, we denote by $\R^E$ the vectorial space of real-valued functions in $E$ whose $dim(\R^E)=m$. For $F\subset E$, the \emph{incidence vector}  of $F$ is defined as follows:
$$\chi_F(e) = \left\{ \begin{array}{cl}
                      1, & \mbox{if} \ e\in F;\\
                      0, & \mbox{otherwise}. \\
                      \end{array}  \right.$$

In general, we identify each subset of edges with its respective incidence vector. The \textit{matching polytope}  of $G$, $\mathcal{M}(G)$, is the convex hull of the incidence vectors of  the matchings in $G$.

In this paper, we are interested in studying the skeleton of a polytope obtained from matchings of a given graph. For more definitions and notations of polytopes, see \cite{Branko}.

\begin{figure}[!ht]
      \centering  
      \includegraphics[scale=0.19]{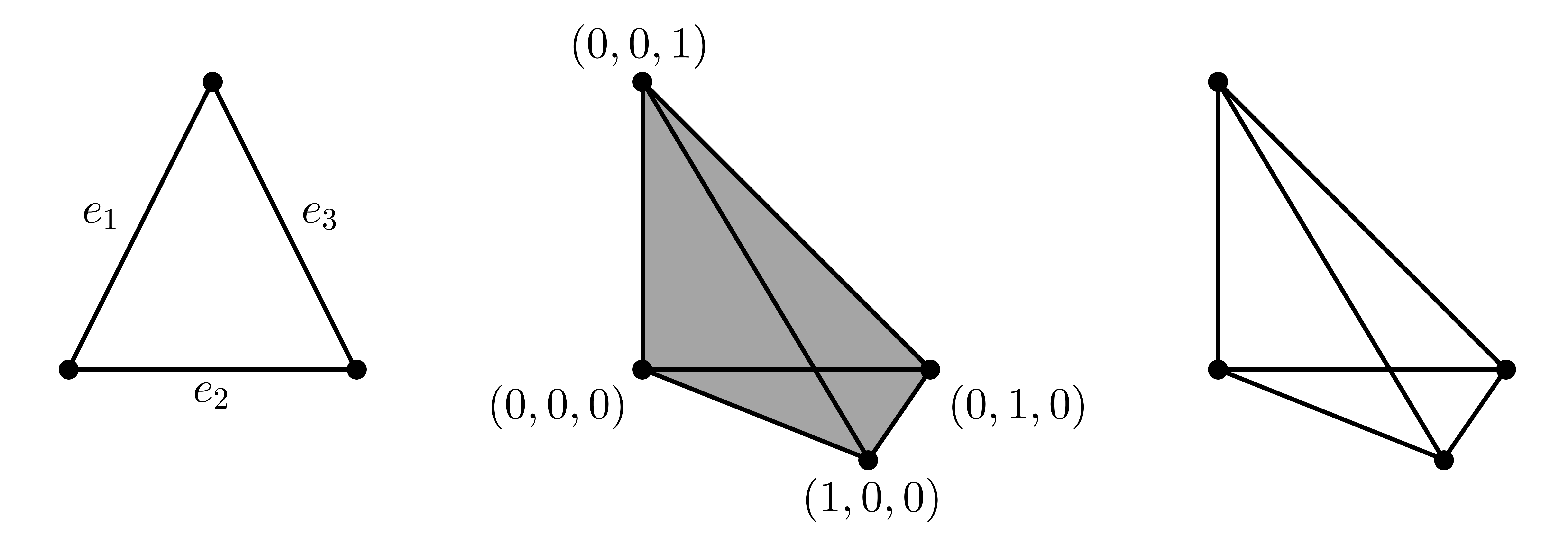}
	\caption{$K_3$ and $\mathcal{G}(\mathcal{M}(K_3)) \equiv K_4$}
	\label{figura57-1}
\end{figure}

The graph $G \equiv K_3$ has the following matchings: $\varnothing$, $\{e_1\}$, $\{e_2\}$ and $\{e_3\}$, where $e_1$, $e_2$ and $e_3$ are edges of $G$. Figure \ref{figura57-1} displays the matching polytope of $K_3$, $\mathcal{M}(K_3)=conv\{(0,0,0), (1,0,0), (0,1,0), (0,0,1)\}$ which corresponds to a tetrahedron in $\R^3$. Its skeleton is $\mathcal{G}(\mathcal{M}(K_3)) \equiv K_4$.

Two matchings $M$ and $N$ are said \textit{adjacent}, $M \sim N$, if and only if their correspondent vertices $\chi_M \equiv M$ and $\chi_N \equiv N$ are adjacent in the skeleton of the matching polytope. The \textit{degree of a matching} $M$,  denoted $d(M)$, is the degree of the correspondent vertex in $\mathcal{G}(\mathcal{M}(G))$. Next theorems characterize the adjacency of two matchings, $M$ and $N$, by their symmetric difference $M \Delta N=(M\backslash N) \cap (N\backslash M)$.

\begin{teo}(\cite{Chvatal}) \label{RP-11}
Let $G$ be a graph. Two distinct matchings $M$ and $N$ of $G$ are adjacent in the matching polytope $\mathcal{M}(G)$ if and only if $M\Delta N$ is a connected subgraph of $G$.
\end{teo}

\begin{teo}(\cite{Schrijver3})\label{CB-14}
Let $G$ be a graph. Two distinct matchings $M$ and $N$ in $G$ are adjacent in the matching polytope $\mathcal{M}(G)$ if and only if $M\Delta N$ is a path or a cycle in $G$.
\end{teo}

Note that the symmetric difference $M \Delta N$ of two adjacent vertices of $\mathcal{M}(G)$ is an $M$ (and $N$)-alternating path (or, cycle). And if it is a cycle, it is an even cycle.

\begin{figure}[!h]
      \centering  
      \includegraphics[scale=0.2]{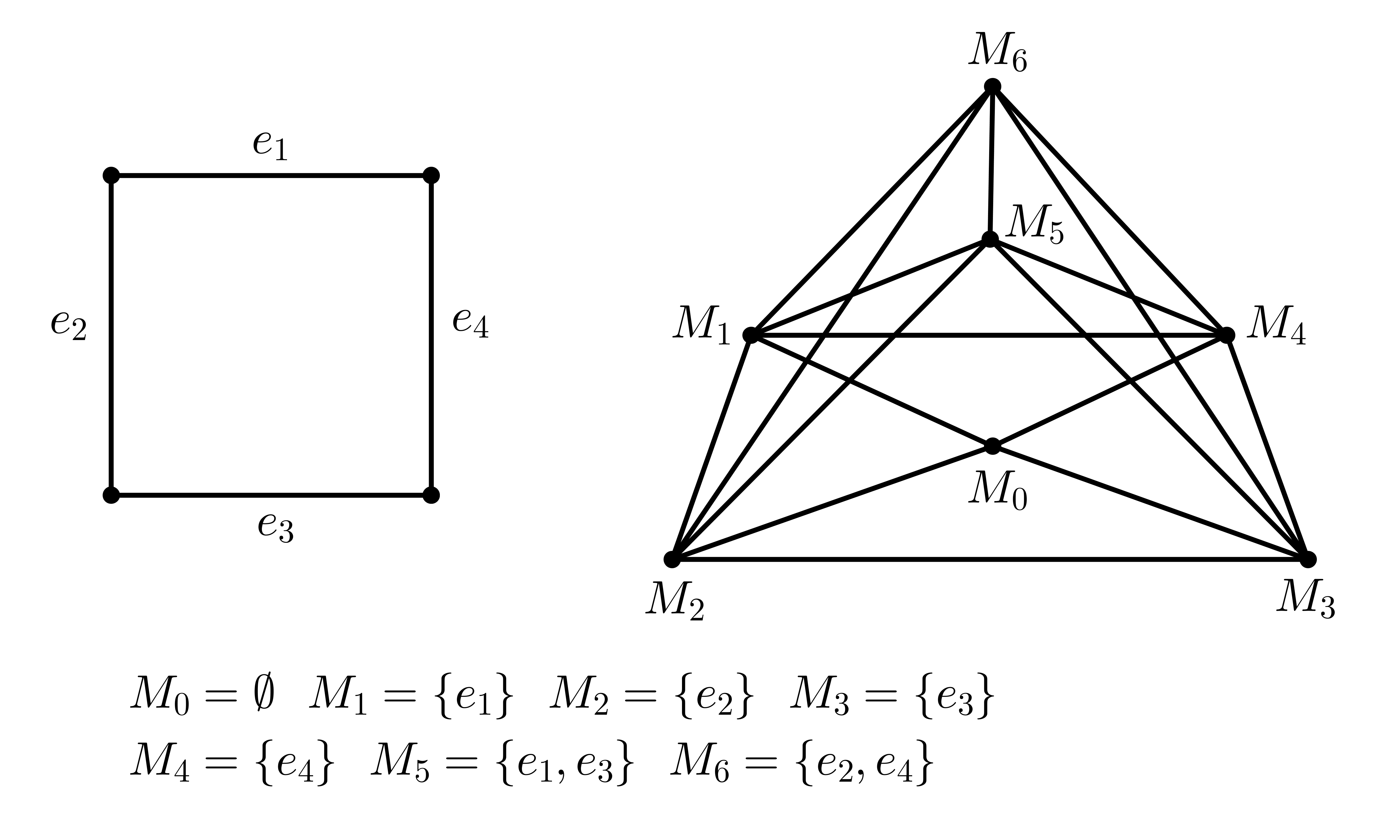}
	\caption{$C_4$ and $\mathcal{G}(\mathcal{M}(C_4))$}
	\vspace{-.19cm}
	\label{figura83}
\end{figure}

Figure \ref{figura83} shows the skeleton of the matching polytope of  $C_4$. Since  $M_1\Delta M_2$ is a path and $M_5\Delta M_6$ is a cycle, $M_1 \sim M_2$ and $M_5 \sim M_6$. However, once $M_1\Delta M_3$ is a disconnected subgraph of $G$, $M_1 \not\sim M_3$.

Let $T$ be a tree with $n$ vertices. The acyclic Birkhoff polytope, $\Omega_n(T)$, is the set of $n \times n$ doubly stochastic matrices $A=[a_{ij}]$ such that the diagonal entries of $A$ correspond to the vertices of $T$ and each positive entry of $A$ is either on the diagonal or on a correspondent position of each edge of  $T$. The matching polytope $\mathcal{M}(T)$ and the acyclic Birkhoff polytope $\Omega_n(T)$  are affinely isomorphic, see \cite{Liliana}. The skeleton of $\Omega_n(T)$ was studied in \cite{Nair, Rosario} where some results about the structure of such graph were presented. In the sequence, we highlight the following contributions given in those papers.

\begin{teo}(\cite{Nair})\label{CB-18}
If $T$ is a tree with $n$ vertices, the minimum degree of $\mathcal{G}(\mathcal{M}(T))$ is $n-1$.
\end{teo}

\begin{teo}(\cite{Nair})\label{pend} Let $T$ be a tree with $n$ vertices and $M$ be a matching of $T$. The degree of $M$ in $\mathcal{G}(\mathcal{M}(T))$ is $d(M)=n-1$ if and only if $M$ is either the empty matching or every edge of $M$ is a pendent edge of $T$. \end{teo}

In this paper, we generalize the results above for an arbitrary graph $G$. In next section, we obtain a formulae to compute the degree of a vertex of $\mathcal{G}(\mathcal{M}(G))$. In Section 3, we prove that the degree of a matching $M$ is non decreasing function of the cardinality of $M$. As a consequence, we conclude that the minimum degree of $\mathcal{G}(\mathcal{M}(G))$ is equal to the number of edges of $G$. Section 4 ends the paper with two theorems of characterization: the first gives a necessary and sufficient condition under $G$ in order to have $\mathcal{G}(\mathcal{M}(G))$ as a regular graph and, the second identifies those matchings of $G$ for which their correspondent vertices of the skeleton have the minimum degree.

\section{The degree of a matching in a graph.}

As we have pointed in the previous section, the first advances on the computation of degree of a matching of tree were due to Abreu et al. \cite{Nair} and Fernandes \cite{Rosario}. In the present section we generalize those results by presenting formulae for the computation of the degree of a vertex of the skeleton of the matching polytope of an arbitrary graph. First, we enunciate two simple results that are immediate consequences of Theorem \ref{CB-14}. In these, we denote the cardinality of a set M by $|M|$.

\begin{prop}\label{0ou1}
Let $M$ and $N$ be matchings of a graph $G$. If $M$ is adjacent to $N$ in $\mathcal{G}(\mathcal{M}(G))$ then $||M|-|N|| \in \{ 0, 1 \}.$
\end{prop}

\begin{pf} Let  $M$ and $N$ be matchings of $G$ such that $M \sim N$ in $\mathcal{G}(\mathcal{M}(G))$.
From Theorem \ref{CB-14}, $M\Delta N$ is a $M$ (and $N$)-alternating path (or, cycle) in $G$. If the length of the path (cycle) is even, then $|M|=|N|$. Otherwise, $||M|-|N|| = 1$. \end{pf}

\begin{prop}\label{RP-12}
Let $G$ be a graph with $m$ edges. The degree of the empty matching is equal to the number of edges of $G$, that is, $d(\varnothing)=m.$
\end{prop}

\begin{pf} Let $G$ be a graph with $m$ edges and suppose $M \neq \varnothing$ be a vertex of $\mathcal{G}(\mathcal{M}(G))$ such that
$M \sim \varnothing$. We know that $M \Delta \varnothing = M$ and, so, from Proposition \ref{0ou1}, $|M|=1$. For some $e\in E(G)$, $M =\{e\}$ is an one-edge matching and, then, $d(\varnothing ) \leq m$. The reciprocal, $m \leq d(\varnothing)$, comes from  Theorem \ref{CB-14},  once the one-matching edge is a path. \end{pf}

Let $G$ be a graph with a matching $M$ and $P$ be an $M$-alternating path with at least two vertices. We say that:\\
(i) $P$ is an \textit{$oo$-$M$-path} if its pendent edges belong to $M$;\\
(ii) $P$ is a \textit{$cc$-$M$-path} if its pendent vertices are both $M$-unsaturated;\\
(iii) $P$ is an \textit{$oc$-$M$-path} if one of its pendent edges belongs to $M$ and one of its pendent vertex is $M$-unsaturated.

An $M$-alternating path $P$ is called an \textit{$M$-good path} if and only if $P$ is one of those paths defined above.

The concept of $M$-good path, introduced by \cite{Nair}, has a important role to determining the degree of a matching of a tree. Since it is our goal to deduce a formulae for the degree of a matching $M$ of any graph, we must consider $M$-alternating cycles. In order to facilitate the writing of the proofs that follow, we will call such cycles by \emph{$M$-good cycles}. Finally, note that an $M$-good cycle $C$ has a perfect matching, given by $M \cap E(C)$.

\begin{figure}[!h]
      \centering  
      \includegraphics[scale=0.9]{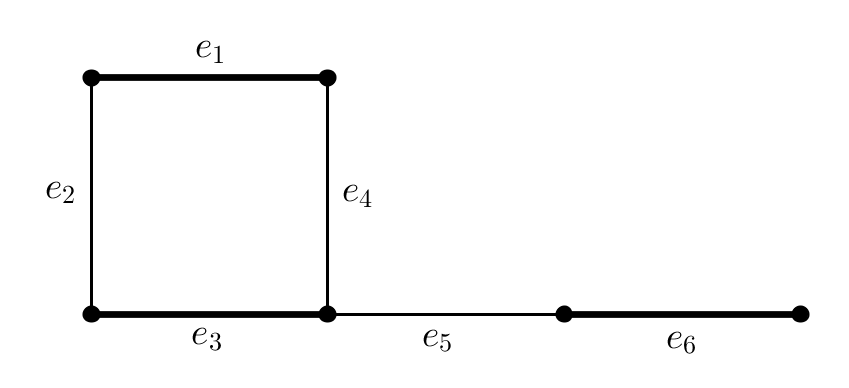}
      \vspace{-.19cm}
	\caption{ $C=e_1e_2e_3e_4$ is $M$-good cycle, for the matching $M=\{e_1,e_3,e_6\}$. }
	\label{figura53-1}
\end{figure}

In Figure \ref{figura53-1}, $M=\{e_1,e_3,e_6\}$ is a perfect matching of the graph. Hence, there are neither $cc$-$M$-paths nor $oc$-$M$-paths in $G$. However, $e_1e_2e_3$ is an $oo$-$M$-path of $G$ but $e_1e_2$ is not an $M$-good path. Moreover, $e_1e_2e_3e_4$ is an $M$-good cycle.

Denote $\nu_{oo}(M)$, $\nu_{cc}(M)$ and $\nu_{oc}(M)$ the numbers of $oo$-$M$-paths, $cc$-$M$-paths and $oc$-$M$-paths, respectively.  Denote $\nu_P(M)$ the number of $M$-good paths of $G$. So, $\nu_P(M) = \nu_{oo}(M) + \nu_{cc}(M)+ \nu_{oc}(M)$. Similarly, the number of $M$-good cycles of $G$ is denoted by $\nu_C(M)$.

\begin{teo}\label{grauM}
Let $M$ be a matching of a graph $G$. The degree of $M$ in the skeleton $\mathcal{G}(\mathcal{M}(G))$ is given by
\begin{equation}\label{eq:grauM}
d(M) = \nu_P(M) + \nu_C(M)~.
\end{equation}
\end{teo}

\begin{pf} Let $M$ and $N$ be adjacent vertices in $\mathcal{G}(\mathcal{M}(G))$. From Theorem \ref{CB-14},  $M\Delta N$ is an $M$-alternating path or a cycle in $G$. If $M \Delta N$ is a path $P$, then it is an $M$-good path. In fact, supposing that $P=e_1e_2\ldots e_{k}$ is not an $M$-good path. Then, two cases can happen. In the first one, $e_1\notin M$ and $e_{k}\notin M$. Besides, one of these edges has an $M$-saturated terminal vertex $u$. So, there is an edge $f \in M\backslash P$ such that $f$ is incident to $u$. Once $P=M\Delta N$ and $f\notin M\Delta N$, then $f\in N$. Consequently,  $N$ is not a matching of $G$. In the second case, suppose $e_1 \in M$ and $e_k \notin M$. Similarly, we get the same contradiction. So, $P$ is an $M$-good path.

Now, assume $P$ as an $M$-good path (or cycle) in $G$ and, do  $N=M\Delta P$. If $P$ is a cycle, it is clear that $N$ is also a matching of $G$. If $P$ is a path, it is also occurs. Indeed, in this case, if $N$ is not a matching of $G$, there is two adjacent edges $f$, $g \in N$ such that $f \in M\backslash P$ and $g \in P\backslash M$. Since $P$ is an $M$-alternating path, $g$ is pendent edge of $P$ and its pendent vertex is $M$-saturated, which leads us to a contradiction. So, $N$ is a matching of $G$. As $M\Delta N=P$ is a path (or cycle) of $G$, from Theorem \ref{CB-14}, $M\sim N$.

Finally, for $N$ and $N'$ matchings of $G$, we have $M\Delta N=M\Delta N'$ if and only if $N=N'$. Hence, $d(M)=\nu_{oo}(M) + \nu_{cc}(M) + \nu_{oc}(M) + \nu_C(M). $\end{pf}

Proposition \ref{RP-12} is also a consequence of Theorem \ref{grauM}. Assume $M=\varnothing$ and $N \sim M$. Then, $M\Delta N$ is one-edge matching and, so, it is a $cc$-$M$-path of $G$.

A matching  $M = \{e_1,  \dots, e_k \} $  is said to be a \emph{matching without common neighbors} of $G$ if and only if
\begin{center}
$\forall i, j \in \{1,\dots, k\},  ~e_i,~e_j \in M,  ~e_i \neq e_j ~\Rightarrow  I(e_i) \cap I(e_j)= ~~ \varnothing.$
\end{center}

\begin{lema}\label{wnccycle}
Let $M$ be a matching without common neighbors in a graph $G$. If $C$ is a cycle of $G$ then $C$ is not an $M$-good cycle.
\end{lema}

\begin{pf} Let $M$ be a matching of $G$ without common neighbors. Suppose there is an $M$-good cycle $C$ of $G$ given by sequence $e_1 f_1e_2f_2 \ldots e_t f_t$, where $f_j \notin M $ and $e_j\in M$, $1 \leq j \leq t$. Since $E(C) \cap M$ is a perfect matching in $C$, there is $i$, $1 \leq i \leq t-1$, such that $e_if_ie_{i+1}$ is a path. Hence, $\{f_i\} \subseteq I(e_i) \cap I(e_{i+1})$.  Once $M$ is a matching without common neighbors, this is an absurd. \end{pf}

\begin{lema}\label{wcnpath}
Let $G$ be a graph and $M$ be a matching without common neighbors. If $P$ is an $M$-good path then $P$ is an $M$-alternating path with length is at most $3$ in which at most one edge belongs to $M$.
\end{lema}

\begin{pf} The proof follows straightforward from $M$ to be a matching without common neighbors and $P$ be an $M$-alternating path of $G$. \end{pf}

\begin{teo}\label{semviz}
Let $M=\{e_i=u_iv_i |  u_i, v_i \in V, 1 \leq i \leq s \}$ be a matching  without common neighbors of a graph $G$. The degree of $M$ is
\begin{equation}\label{eq:grauM2}
d(M)=k+\sum_{i=1}^{s} \left(d(u_i)d(v_i)-|N(u_i)\cap N(v_i)|\right),
\end{equation}
where $k\geq 0$ is the number of edges of $G$ which neither incides to $u_i$ nor to $v_j$, \ for all $1\leq i < j \leq s$.
\end{teo}

\begin{pf} Let $M=\{e_i=u_iv_i |  u_i, v_i \in V, 1 \leq i \leq s \}$ be a matching without common neighbors of $G$. From Lemma \ref{wnccycle}, $G$ does not have $M$-good cycles and, from Lemma \ref{wcnpath}, if $P$ is an $M$-good path, $P$ has at most length $3$ with at most one edge of $M$. Hence, if $P$ is an $M$-good path of $G$, $P$ has to take one of the cases bellow.
\begin{description}
\item (1) $P$ is an $oo$-M-path. Then, $P=e$ such that $e\in M$. There are $s$ of these paths in $G$;
\item (2) $P$ is an $oc$-M-path. So, $P=ef$ such that $e\in M$ and $f\notin M$. Of course, for some $i$, $1\leq i\leq s$ we have $e=u_iv_i$ and $f$ is incident to $u_i$ or to $v_i$. There are $\left(d(u_i)-1\right)+\left(d(v_i)-1\right)$ of these paths;
\item (3) If $P$ is a $cc$-M-path, we have to consider two possibilities for $P$.  Firstly, $P=feg$ with $e \in M$ and $f, g \notin M $. So, for some $i$, $1\leq i \leq s$, $e=u_iv_i$ such that $f$ is incident to $u_i$ and $g$ is incident to $v_i$. In this case, there are $\left(d(u_i)-1\right)\cdot \left(d(v_i)-1\right)-|N(u_i)\cap N(v_i)|$ of such paths.
The second possibility is $P=f$ with $f \notin M$. Since $P$ is a $cc$-M-path, for each $e \in E(G)$ that is incident to $f$, $e \notin M$. We can admit that there are $k$ of these edges in the graph that satisfy this last case.
\end{description}
From the itens (1), (2) and (3) and, by applying Theorem \ref{grauM}, we obtain
$$d(M)= s + k + \sum_{i=1}^{s} ((d(u_i)-1)+(d(v_i)-1) + (d(u_i)-1)(d(v_i)-1)-|N(u_i)\cap N(v_i)|)=$$
$$= k+\sum_{i=1}^{s} (d(u_i)d(v_i)-|N(u_i)\cap N(v_i)|).$$ \end{pf}

\section{Vertices with minimum degree in the skeleton.}

We begin this section proving that, if a matching within the other, the degree of the first is at most equal to the degree of the last. Based on this, we prove that the minimum degree of $\mathcal{G}(\mathcal{M}(G))$ is equal to the number of edges of $G$. The section follows by determining the degree of a matching of a graph whose components are stars or triangles.

\begin{teo}\label{RP-4}
Let $M$ and $N$ be matchings of a graph $G$. If $M\subset N$, then $d(M)\leq d(N)$.
\end{teo}

\begin{pf} Let $M$ and $N$ be matchings of a graph $G$ such that $N=M \cup \{e\}$, where $e \in E(G)$. Consider $\mathbb{B}_M$ the sets of $M$-good paths and $M$-good cycles of $G$. Similarly, define  $\mathbb{B}_N$. From Theorem \ref{grauM}, $d(M)=|\mathbb{B}_M|$ and $d(N)=|\mathbb{B}_M|$.  Build the function $\varphi_e: \mathbb{B}_M \longrightarrow \mathbb{B}_N$ such that, for every cycle $C \in \mathbb{B}_M$, $\varphi_e(C) = C$ and, for every path $P \in \mathbb{B}_M,$

\begin{center}
$\varphi_e(P) = \left \{ \begin{array}{lll}
                                         P, & \mbox{if $e \in P$;} \\
                                         P, & \mbox{if $e \notin P$ and if $f \in P$ then not $f \sim e$;} \\
				         P \cup \{ e\}, & \mbox{otherwise.}
				       \end{array}
				       \right.$
				       \end{center}
				       				       				
By construction, for distinct paths or cycles belonging to $\mathbb{B}_M$, we have distinct images in $\mathbb{B}_N$. Then, $\varphi_e$ is a injective function and so,  $d(M) \leq d(N)$.

In the general case, do $N \backslash M = \{e_1, e_2,\ldots, e_k\}$, $N_1 =M  \cup \{e_1\}, N_2= M \cup \{e_1,e_2\},\ldots,$  and $N_k=N$. By the same argument used before, we get  $d(M) \leq d(N_1), d(N_1) \leq d(N_2), \ldots,$ and $d(N_{k-1}) \leq d(N)$. Consequently, $d(M) \leq d(N).$ \end{pf}

\begin{teo}\label{RP-15}
Let $G$ be a graph with $m$ edges. The minimum degree of $\mathcal{G}(\mathcal{M}(G))$ is equal to $m$.
\end{teo}
\begin{pf} Let $M$ be a matching of a graph $G$.  Since $\varnothing \subseteq M$, from Theorem \ref{RP-4}, $d(\varnothing) \leq  d(M)$. Besides, by Proposition \ref{RP-12}, $d(\varnothing) =m$ and, the result follows. \end{pf}

 An edge $e=uv$ of a graph $G$ is called a \textit{bond} if $d(u)=d(v)=2$ and $|N(u)\cap N(v)|=1$. Note that if $e$ is a bond of $G$, $e$ is an edge of a triangle of graph. However, the reciprocal is not necessarily true.

The following lemma allows us to obtain, in the next section, a characterization of graphs (connected or not) for which the respective skeletons are regular. See that the skeleton of a polytope is always a connected graph even if the original graph is disconnected, \cite{Goodman}.

\begin{lema}\label{lema1}
Let $G$ be a graph with $m$ edges such that $G$ is a disjoint union of stars and triangles. For every $M$, a matching of $G$,  $d(M)=m$.
\end{lema}

\begin{pf} Let $G$ be a graph with $m$ edges. Suppose $M$ a matching of $G$. If $M = \varnothing$,  the result follows straightforward from Proposition \ref{RP-12}.

From hypotheses, $G$ is union of $r$ triangles $K_3$ and $p$ stars $S_{1,t_j}, 1\leq j \leq p$. So, it can be written as

\begin{equation}\label{eq:sum}
G=rK_3 \bigcup\limits_{j=1}^{p} S_{1,t_j},
\end{equation}
for some non negative integers $r$, $p$ and $t_j, 1\leq j \leq p$  such that

\begin{equation}\label{eq:edges}
m=3r+\sum_{j=1}^{p} t_j.
\end{equation}

Since (\ref{eq:sum}) holds, each two distinct edges of $M$ belong to distinct components of $G$. Also, $M$ is a matching without common neighbors of $G$. Hence, the expression (\ref{eq:sum}) can be rewritten as

\begin{equation}\label{eq:sum2}
G = r_1K_3 \bigcup\limits_{j=1}^{p_1} S_{1,t_j} \cup r_2K_3 \bigcup\limits_{z=1}^{p_2} S_{1,t_z},
\end{equation}

where $r_1+p_1$ is the cardinality of $M$ and $r_2+p_2$ is the number of components without any edges of $M$. So, the number of edges that not are incident to any edge of M is
\begin{equation}\label{sete}
k=3r_2+ \sum_{z=1}^{p_2} t_z.
\end{equation}

For each $i, 1 \leq i \leq r_1+p_1$, let $e_i=u_iv_i \in M$, denote $s_i = d(u_i)d(v_i) -|N(u_i) \cap N(v_i)|$. If $e_i \in E(K_3)$ then $s_i=3$. Otherwise, $e_i \in E(S_{1,t_j})$ is a pendent edge of $G$ and $s_i = t_j$.

As $M$ is a matching without common neighbors in $G$,  from (\ref{sete}) and applying Theorem \ref{semviz}, we get
\begin{equation}\label{oito}
d(M)=(3r_2+\sum_{z=1}^{p_2} t_z) + \sum_{i=1}^{r_1+p_1} s_i = 3(r_1+r_2)+ \sum_{j=1}^{p_1} t_j + \sum_{z=1}^{p_2} t_z .
\end{equation}

Once $r_1+r_2=r$ and $p_1+p_2=p$, from (\ref{eq:edges}), we get $d(M)=m$.  \end{pf}

\vspace{0.2cm}

\begin{figure}[!h]
      \centering  
      \includegraphics[scale=.95]{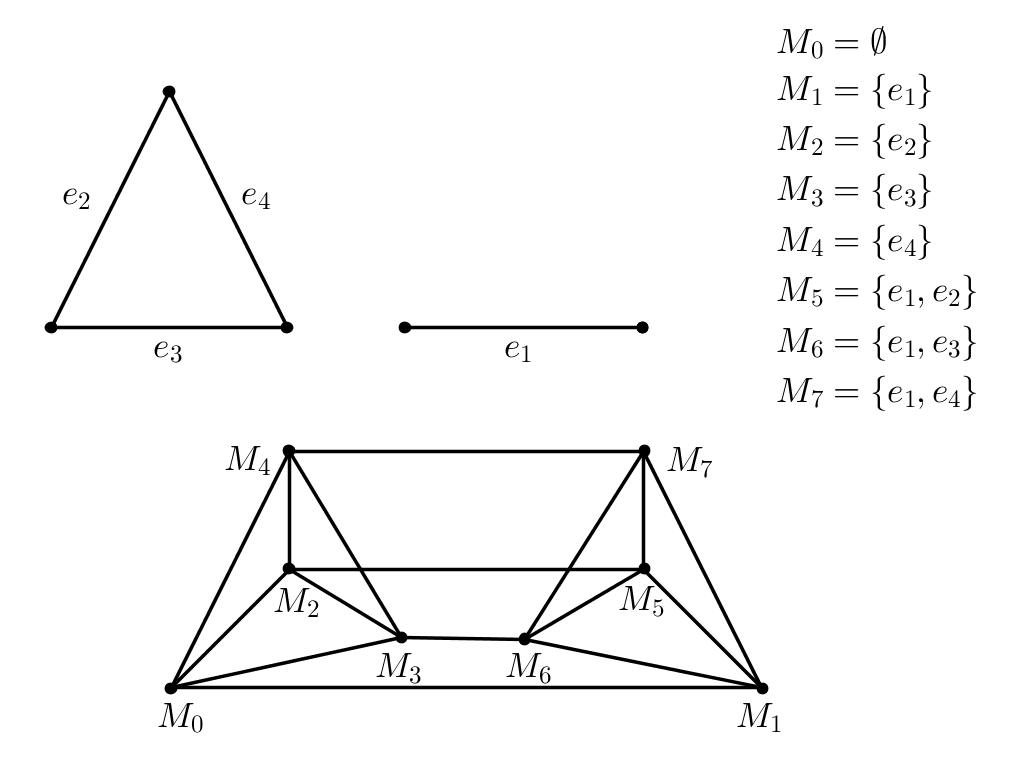}
      \vspace{-.19cm}
	\caption{$G=K_3\cup S_{1,1}$ and $\mathcal{G}(\mathcal{M}(G))$}
	\label{figura81}
\end{figure}

 Figure \ref{figura81} displays the graph $G=K_3\cup S_{1,1}$ and its skeleton, $\mathcal{G}(\mathcal{M}(K_3\cup S_{1,1}))$.
 For $M = \{e_1, e_2\}$,  we have $r =r_1=1$ and $r_2 =0$, $p =p_1= 1$ and $p_2 =0$. Besides, $k=0$, $s_1 = 1$ and $s_2 = 3$. So, $d(M)=k+s_1+s_2=0+1+3=4$.

\begin{prop}\label{minimal}
Let $G$ be a graph with $m$ edges and $M=\{e\}$ be an one-edge matching of $G$.  The degree of $M$ is $d(M) = m$ if and only if $e$ is either a bond or a pendent edge of $G$.
\end{prop}

\begin{pf} Let $G$ be a graph with $m$ edges and $e=uv$ an edge of $G$. We know that $m - d(u) - d(v) + 1$ is the number of edges that neither is incident to $u$ nor to $v$. From Theorem \ref{semviz}, $d(\{e\})=m - d(u) - d(v) + 1 + d(u)d(v) - t$, where $t = |N(u) \cap N(v)|$. But, $d(\{e\}) = m$ if and only if $d(u)d(v) -  t = d(u) + d(v)  - 1$, i.e., $(d(u) - 1)(d(v) - 1) = t$. Since $d(u) > t$ and  $d(v) > t$, $(d(u) -1)(d(v) -1)  \geq t^2$. Moreover, $(d(u) -1)(d(v)- 1) = t$ if and only if $t = 0$ or $t = 1$. In the first case, $e$ is a pendent edge and, otherwise, $d(u)=d(v)=2$ and so, $e$ is a bond. \end{pf}

Note that, if $G$ is a graph without bonds and pendent edges, $\mathcal{G}(\mathcal{M}(G))$ has only one vertex $M=\varnothing$ with the minimum degree. It is not
difficult to see that all 2-connected graphs satisfy this property.


\section{Characterizing regular skeletons and matchings with minimum degree}

We close this paper with two characterizations. The first gives a necessary and sufficient condition under a graph in order to have the skeleton of its matching polytope as a regular graph. The second identifies all matchings of a graph which correspondent vertices in the skeleton have the minimum degree.

\begin{prop} A graph $G$ with $m$ edges is a disjoint union of stars and triangles if and only if $\mathcal{G}(\mathcal{M}(G))$ is an $m$-regular graph.
\end{prop}

\begin{pf} From Lemma \ref{lema1}, if $G$ is a disjoint union of stars or triangles, the skeleton of its matching polytope is a regular graph. Suppose $\mathcal{G}(\mathcal{M}(G))$ an $m$-regular graph. Then, for every $e \in E(G)$, we have $d(\{e\}) = d(\varnothing)$ and, from Proposition \ref{RP-12},  $d(\{e\}) = m.$ Besides, by Proposition \ref{minimal}, this occurs only if $e$ is a bond or a pendent edge of $G$. \end{pf}

\begin{teo} Let $M \neq \varnothing$ is a matching of a graph $G$ with $m$ edges. The degree of $M$ is $m$, that is, $d(M)=m$ if and only if every edge of $M$ is a bond or a pendent edge of $G$.
\end{teo}

\begin{pf} Let $G$ be a graph with a matching $M$. Suppose there is $e \in M$ such that $e$ is neither a bond nor a pendent edge of $G$. So, from Proposition \ref{minimal} and Theorem \ref{RP-4}, $m < d(\{e\} ) \leq d(M)$. Consequently, $m \neq d(M)$. By the contrapositive, if $d(M) = m$, every edge of $M$ is a pendent edge or a bond of the graph.

Suppose now that $M$ is a matching of $G$ with $s$ edges such that if $e \in M$, $e$ is a bond or $e$ is a pendent edge of $G$. Let $N$ be a matching such that $N \sim M$ in $\mathcal{G}(\mathcal{M}(G))$. From here and by Theorem \ref{CB-14}, $M \Delta N$ is an $M$-good path $P$ or an $M$-good cycle $C$ of $G$. Since $C$ is an even cycle, $C \neq K_3$. So, $C$ does not have bonds. Consequently, $M \Delta N$ is a path.

Concerning the path $P$, only their pendent edges can belong to $M$. Moreover, once $P$ is an alternated path, it has length at most length 3. Hence, there
are only the following possibilities to $P$:

\begin{enumerate}

\item If $P$ is an $oo$-M-path, then $P = e$, where $e \in M$, or $P = e_1fe_2$, where $f  \notin M$ and $I(f) \cap M = \{e_1, e_2\}$. In the first case, there are $s$ possibilities to $P$ and, in the second, there are $t_1$ possibilities to $P$, where $t_1$ is the number of edges of $E(G) \backslash M$ such that the both terminal vertices are incident to an edge of $M$;

\item If $P$ is a $cc$-M-path, $P = f$, where $f \notin M$ and $I(f) \cap M = \varnothing$. Here, there are $t_2$ possibilities to $P,$ where $t_2$ is the number of edges of $E(G) \backslash M$ for which every edge of $M$ does not incident to the any end vertices of those edges;

\item Finally, if $P$ is an $oc$-M-path, then $P = ef$, where $f \notin M$ and $I(f) \cap M = \{e\}$. In this last case, there are $t_3$ possibilities to $P$, where $t_3$ is the number of edges of $E(G) \backslash M$ which only one end vertex of edge is incident to some edge of $M$.

\end{enumerate}

From (1), (2) and (3) possibilities above, $s +\sum\limits_{i=1}^3 t_i = s + |E(G) \backslash M| = m$ is the number of the possibilities to have $P$ as an M-good path of $G$. From Theorem \ref{grauM} it follows that $d(M) = m$. \end{pf}

\stepcounter{section}

\end{document}